\documentclass[12pt]{amsart}
\usepackage{amssymb}

\begin{document}

\numberwithin{equation}{section}

\theoremstyle{plain}
\newtheorem{thm}[equation]{Theorem}
\newtheorem{lemma}[equation]{Lemma}
\newtheorem{cor}[equation]{Corollary}
\newtheorem{prop}[equation]{Proposition}
\newtheorem{conj}[equation]{Conjecture}

\theoremstyle{definition}
\newtheorem{defn}[equation]{Definition}
\newtheorem{notation}[equation]{Notation}

\theoremstyle{remark}
\newtheorem{rem}[equation]{Remark}
\newtheorem{exam}[equation]{Example}
\newtheorem{ack}[equation]{Acknowledgments}

\newenvironment{thmref}{\thmrefer}{}
\newcommand{\thmrefer}[1]{\renewcommand\theequation
  {\protect\ref{#1}$'$}\addtocounter{equation}{-1}}

\newcommand{\Cay}{\overrightarrow{\operatorname{Cay}}}
\newcommand{\Z}{\mathbb{Z}}
\newcommand{\cyc}{\overrightarrow{C}}

\newenvironment{basecase}[1][\unskip]{
  \medskip \noindent \em Base case.\ }{\unskip\upshape}

\newenvironment{inductionstep}[1][\unskip]{
  \medskip \noindent \em Induction Step.\ }{\unskip\upshape}

\title{Hamiltonian Paths in Cartesian Powers of Directed Cycles}

\author{David Austin}

\address{
Department of Mathematics and Statistics,
Grand Valley State University,
Allendale, MI  49401}

\author{Heather Gavlas}

\address{
Department of Mathematics and Statistics,
University of Vermont,
Burlington, VT 05401}

\author{Dave Witte}

\address{
Department of Mathematics,
Oklahoma State University,
Stillwater, OK  74078
}
\maketitle

\begin{abstract}
  The vertex set of the $k^{\text{th}}$ cartesian power of a directed 
cycle of length~$m$
can be naturally identified with the abelian group $(\Z_m)^k$. For 
any two elements $v =
(v_1,\ldots,v_k)$ and $w = (w_1,\ldots,w_k)$ of $(\Z_m)^k$, it is 
easy to see that if there
is a hamiltonian path from~$v$ to~$w$, then
  $$ v_1 + \cdots + v_k \equiv w_1 + \cdots + w_k + 1 \pmod{m} .$$
  We prove the converse, unless $k = 2$ and $m$~is odd.
  \end{abstract}

\section{Introduction}

The cartesian product of any number 
of (undirected) cycles always contains a hamiltonian cycle.
(See (\ref{cart_prod_defn}) for the definition of the cartesian 
product.) Work of Chen and
Quimpo \cite{CQ1} implies the following stronger result, which 
provides a simple
characterization of the pairs of vertices that can be joined by a 
hamiltonian path.

\begin{thm}[Chen-Quimpo \cite{CQ1}]
  Let $X$ be the cartesian product of $k$ cycles of lengths 
$m_1,m_2,\ldots,m_k$, with $k \ge 2$
and  each $m_i \ge 3$.
  \begin{enumerate}
  \item If some $m_i$ is odd, then $X$ is hamiltonian connected. That 
is, for any two
vertices $u$ and~$v$ of~$X$, there is a hamiltonian path from~$u$ to~$v$.
  \item If each $m_i$ is even, then $X$ is hamiltonian laceable. That 
is, $X$ is bipartite
and, for any two
vertices $u$ and~$v$ of~$X$, either there is a hamiltonian path 
from~$u$ to~$v$, or there
is a path of even length from~$u$ to~$v$.
  \end{enumerate}
  \end{thm}

It would be interesting to have a similar result in the directed 
case. The first step,
which has been completed, is to determine which cartesian products of 
directed cycles have
hamiltonian cycles. Rankin \cite{R1} implicitly gave a necessary and 
sufficient condition
for the  existence of a hamiltonian cycle in the cartesian product of 
two directed cycles;
however, this result went unnoticed by graph theorists.  Thirty years 
later, Trotter and
Erd\"{o}s \cite{TE1}  rediscovered the characterization. Curran and 
Witte \cite{CW1} showed
that there is a hamiltonian cycle in the cartesian product of three 
or more nontrivial
directed cycles.

\begin{thm}[Rankin \cite{R1}, Trotter-Erd\"{o}s \cite{TE1}, and 
Curran-Witte \cite{CW1}]
\label{ham_cyc}
  Let $X$ be the cartesian product of $k$ directed cycles of lengths 
$m_1, m_2, \ldots,m_k$, with $k
\ge 1$ and  each $m_i \ge 2$. Then there is a 
hamiltonian cycle
in~$X$ if and only if either $k \neq 2$ or there exists a pair of 
relatively prime positive
integers $s_1$  and $s_2$ with $s_1 m_1 + s_2 m_2 = m_1 m_2$.
  \end{thm}

For the case $k = 2$, Curran \cite{CW1} strengthened 
Theorem~\ref{ham_cyc} to obtain a
description of the pairs of vertices that can be joined by a 
hamiltonian path in terms of
the geometric configuration of the lattice points in the plane 
triangle with vertices
$(m_1,0)$, $(0,m_2)$, and $(0,0)$. The result is quite technical so, instead of
stating it here, let us
mention that when $k=2$, there always exist pairs of vertices 
that cannot be joined by
a hamiltonian path. In fact, Curran showed, for each vertex~$v$, that 
no more than half of
the vertices of~$X$ are the terminal vertex of some hamiltonian path 
starting at~$v$.

In contrast, we conjecture that there
is a hamiltonian path from~$u$ to~$v$ for any two distinct vertices 
$u$ and~$v$ if 
$\gcd(m_1, m_2, \ldots,m_k) = 1$ and $k \ge 3$. This is a
special case of the following conjecture, which has no restriction on 
$m_1, m_2, \ldots, m_k$.

\begin{notation}
  For vertices $u$ and~$v$ in a (strongly connected) digraph~$X$,  
let $d_X(u,v)$  denote
the length of the shortest directed path from~$u$ to~$v$.
  \end{notation}

\begin{conj} \label{conjecture}
   Let $X$ be the cartesian product of $k$ directed cycles of lengths 
$m_1, m_2, \ldots,m_k$, with $k
\ge 3$ and each $m_i \ge 2$. For vertices $u$ and $v$ 
of  $X$, there is a
hamiltonian path from~$u$ to~$v$ if and only if
  \begin{equation} \label{congruence}
  d_X(u,v) \equiv -1 \pmod{\gcd(m_1, m_2, \ldots,m_k)} .
  \end{equation}
  \end{conj}

  For $X$ as in Conjecture~\ref{conjecture}, the lengths of any two 
directed paths from~$u$
to~$v$ are congruent modulo $\gcd(m_1, m_2, \ldots,m_k)$. This elementary 
observation implies
that \eqref{congruence} is a necessary condition for the existence of 
a hamiltonian path
from~$u$ to~$v$. Our conjecture is that it is also sufficient.

In this paper, we prove Conjecture \ref{conjecture} in the special case where all 
of the directed cycles
have the same length.

\begin{thmref}{main}
  \begin{thm}
  Let $X$ be the $k^{\text{th}}$ cartesian power of a directed cycle 
of length~$m$, with $k \ge 3$ and $m \ge 2$. For vertices $u$ and $v$ of $X$, there is a 
hamiltonian path 
from~$u$ to~$v$ if and
only if
  $$ d_X(u,v) \equiv -1 \pmod{m} .$$
  \end{thm}
  \end{thmref}

\section{Definitions and Preliminaries}

Let us begin this section with a few definitions.  For background information about the 
graph theoretic terms used below, the reader is directed to \cite{CL1}.  For background 
information on hamiltonian cycles in Cayley digraphs, including the arc-forcing subgroup, 
see the surveys \cite{CG1, WG1}.

\begin{notation}
  The directed cycle on $m$ ($m \ge 2$) vertices is denoted by
$\cyc_m$.
\end{notation}

\begin{defn}\label{cart_prod_defn}
  The \emph{cartesian product} $G = G_1 \times G_2$ of two digraphs $G_1$
and  $G_2$ is the digraph whose vertex set is $V(G) = V(G_1) \times V(G_2)$
and has an arc from  $(u_1, u_2)$ to $(v_1, v_2)$ if and only if either
  $$ \mbox{$u_1 = v_1$ and there is an arc from $u_2$ to $v_2$ in $G_2$} $$
  or
  $$\mbox{$u_2 = v_2$ and there is an arc from $u_1$ to $v_1$ in $G_1$.}$$
\end{defn}

A convenient way of drawing $G_1 \times G_2$ is to first place a copy 
of $G_2$ at each
vertex of $G_1$ and then join corresponding vertices of $G_2$ in the 
copies of $G_2$
placed at adjacent vertices of $G_1$.

\begin{notation}
  In the abelian group $\Z_{m_1} \times  \Z_{m_2} \times \cdots 
\times  \Z_{m_k}$, let
  \begin{align*}
  x_1 &= (1, 0, 0, \ldots, 0) , \\
  x_2 &=  (0, 1, 0, \ldots, 0) , \\
  &\ldots  \\
  x_k &= (0, 0, \ldots, 0, 1) .
  \end{align*}
  The set $S = \{x_1, x_2,  \ldots,  x_k\}$ will denote the \emph{standard 
generating set} of
$\Z_{m_1} \times  \Z_{m_2}  \times  \cdots  \times \Z_{m_k}$.
  \end{notation}

\begin{defn}
  Let $S$ generate a finite (abelian) group $\Gamma$. The \emph{Cayley digraph
$\Cay(\Gamma; S)$} is the digraph whose vertex set is $V(G)=\Gamma$
and has an arc from $g$ to $g+s$ whenever $g \in \Gamma$ and $s 
\in S \smallsetminus
\{0\}$. (We delete $0$ from~$S$ to avoid having loops in the digraph.)

We may write $\Cay(\Z_{m_1} \times  \Z_{m_2}  \times  \cdots  \times 
\Z_{m_k})$, omitting
$S$ from the notation, when the standard generating set is to be used.
  \end{defn}

  \begin{notation}
  A path $P$ in a digraph can be specified by giving an ordered list 
$v_0,v_1,\ldots,v_n$
of the vertices encountered.

In a Cayley digraph  $\Cay(\Gamma;S)$, it is usually more convenient 
to specify the path $P: v_0,v_1,\ldots,v_n$
by
giving the initial vertex~$v_0$ and an ordered list 
$(a_1,a_2,\ldots,a_n)$ of 
elements of~$S$ that label the arcs in $P$ as $v_i = v_0 + 
a_1 + \cdots + a_i$.

  We will also use the notation $(a_1, a_2, \ldots, a_n)^k$  to 
indicate the concatenation
of $k$ copies of the sequence $(a_1, a_2, \ldots, a_n)$.  Thus, for 
example, $(a^2,  b)^3 =
(a, a, b, a, a, b, a, a, b)$.
  \end{notation}

The following observation is well known (and easy to prove).

\begin{prop}
  If $\Gamma = \Z_{m_1} \times  \Z_{m_2} \times \cdots  \times 
\Z_{m_k}$ $(m_i \ge 2)$ and
$S$ is the standard generating set, then $\Cay(\Gamma; S)$ is
isomorphic to the cartesian product of $k$~directed cycles of 
lengths $m_1,$ $m_2,$ $\ldots,$
$m_k$; that is,
  $$\Cay(\Z_{m_1} \times  \Z_{m_2} \times \cdots  \times 
\Z_{m_k}) \cong  \cyc_{m_1} \times \cyc_{m_2} \times
\cdots \times \cyc_{m_k} .$$
  \end{prop}

  Thus, in order to understand cartesian  products of directed cycles, 
it suffices to
understand certain Cayley digraphs. This change of perspective 
provides an algebraic
setting that makes some constructions more transparent. As Cayley digraphs are
vertex-transitive, there is usually no harm in assuming that 
a hamiltonian path starts at the
identity element. In $\Cay 
\bigl( (\Z_m)^k \bigr)$,
there is a hamiltonian path from~$v$ to~$w$ if and only if there is a 
hamiltonian path
from~$0=(0,0,\ldots,0)$ to~$w - v$; thus throughout the rest of this paper, we will 
consider only that hamiltonian paths start at $0$.

The following subgroup of $(\Z_m)^k$, called the ``arc-forcing 
subgroup," is a basic tool
in the study of hamiltonian paths.

\begin{defn}
  For the standard generating set~$S$ of $(\Z_m)^k$, let
  $S^{-} = \{\, -s \mid s \in S\}$ and let
  \begin{align*}
  H &= \langle S+S^{-} \rangle \\
   &= \{\, (v_1,\ldots,v_k) \in (\Z_m)^k \mid v_1 + \cdots + v_k 
\equiv 0 \pmod{m} \,\} .
  \end{align*}
  Then $H$ is  called the \emph{arc-forcing subgroup}.
  \end{defn}

\begin{rem} Let us recall some basic facts.
  \begin{itemize}
  \item For any $s_1,s_2,x \in S$, we have
  $$s_1 - s_2 = (s_1 - x) - (s_2 - x) \in \langle S - x \rangle ,$$
  where $S - x = \{\, s - x \mid s \in S \,\}$. Therefore
  $S - x$ is a generating set for~$H$.
  \item For $v = (v_1, v_2, \ldots, v_k) \in (\Z_m)^k$, we have $v \in 
H - x_1$ if and only
if
  $v_1 + v_2 + \ldots +  v_k \equiv -1 \pmod{m}$.
  \item If there is a hamiltonian path from~$0$ to~$v$ in $\Cay \bigl(
(\Z_m)^k \bigr)$, then $v \in H-x_1$.
  \end{itemize}
  \end{rem}

\section{Hamiltonian Paths}

\begin{notation}
  Throughout this
section:
  \begin{itemize}
  \item $m $ and $k$ are positive integers with $m , k \ge 2$,
  \item $S$ is the
standard generating  set of $(\Z_m)^k$, say $S =\{ x_1, x_2, \ldots, 
x_k\}$, and
  \item $H =
\langle S+S^{-} \rangle$ is the arc-forcing subgroup.
  \end{itemize}
  \end{notation}

In this section, we will prove our main result.

\begin{thm}\label{main}
  If $k \ge 3$, then for every $v \in H - x_1$, there is a 
hamiltonian path from $0$ to $v$
in $\Cay \bigl( (\Z_m)^k \bigr)$.
  \end{thm}

Since $\Cay(H; S - x_1) \cong \Cay \bigl( (\Z_m)^{k-1} \bigr)$, it is 
easy to see that
$\Cay(H; S - x_1)$ contains a  hamiltonian cycle $C$.  If $v$ is any 
element of~$H$ that is
an even distance from $0$ along $C$, then the following proposition 
shows that there is
a hamiltonian path from $0$ to $v - x_1$ in  $\Cay \bigl( (\Z_m)^k 
\bigr)$. This observation
is the main tool in our proof of Theorem~\ref{main}.

\begin{prop}\label{ham_path}
  Let $v \in H$. If there is a hamiltonian cycle~$C$ in $\Cay(H; S - 
x_1)$ such that
$d_C(0,v)$ is even, then $\Cay \bigl( (\Z_m)^k \bigr)$ contains a 
hamiltonian path from $0$
to $v -  x_1$.
  \end{prop}

\begin{proof}
  Let $C : c_0,c_1,\ldots,c_{m^{k-1}}$ be a hamiltonian cycle in $\Cay(H; S - 
x_1)$ with $c_0 = c_{m^{k-1}} = 0$, and define
  $$ \mbox{$\phi \colon \Z_m \times \Z_{m^{k-1}} \to (\Z_m)^k$ by 
$\phi(i,j) = i x_1 +
c_j$.} $$
  For convenience, let $a = (1,0)$ and $b = (1,1)$ in $\Z_m \times 
\Z_{m^{k-1}}$. Then, for
any $v = (i,j) \in \Z_m \times \Z_{m^{k-1}}$, we have
  $$ \phi(v+a) - \phi(v) = x_1 \in S$$
  and
  $$ \phi(v+b) - \phi(v) = x_1 + (c_{j+1} - c_j) \in S ,$$
  as $c_{j+1} - c_j \in S - x_1$. So $\Cay \bigl( (\Z_m)^k \bigr)$ 
contains arcs from $\phi(v)$
to $\phi(v+a)$ and $\phi(v+b)$. Therefore, $\phi$ embeds
  $\Cay \bigl( \Z_m \times \Z_{m^{k-1}}; \{a,b\} \bigr)$
  as a spanning subdigraph of $\Cay \bigl( (\Z_m)^k) \bigr)$.
  Now, for $0 \le n < m^{k-1}/2$, the path
  $$\bigl( (a^{m-2}, b^2)^n, (a^{m-2}, b, a)^{m^{k-1} - 2n - 1}, 
(a^{m-2}, b^2)^n, a^{m-2},
b \bigr) $$
  is a hamiltonian path~$P$ from~$0$ to~$(-1,2n)$ in
  $\Cay \bigl( \Z_m \times \Z_{m^{k-1}}; \{a,b\} \bigr)$.
  So $\phi(P)$ is a hamiltonian path from $\phi(0) = 0$ to
  $\phi(-1,2n) = c_{2n} - x_1$ in $\Cay \bigl( (\Z_m)^k \bigr)$.
  \end{proof}

The following corollary establishes Theorem~\ref{main} in the case 
where $m$ is even.

\begin{cor} \label{main_even}
  If $m$ is even, then for every $v \in H - x_1$, there is a 
hamiltonian path from~$0$
to~$v$ in $\Cay \bigl( (\Z_m)^k \bigr)$.
  \end{cor}

\begin{proof}
  Let $v = (v_1, v_2, \ldots, v_k)$.  Since $v \in H - x_1$,  we have 
that $v_1 + v_2 +
\ldots + v_k \equiv -1 \pmod{m}$.  Since $m$ is even, it  must be the 
case that $v_i$ is odd
for some $i$; without loss of generality, say $v_1$ is odd.  We may assume this without 
loss of generality as $S-x_i$ is a generating set for $H$ and $v \in H-x_i$ for each $i\ (1 
\le i \le k)$. Let $C$ 
be a hamiltonian
cycle  in  $\Cay(H;S-x_1)$.

For  $\epsilon \in \{0, 1\}$, let
  $$H_\epsilon = \{ h \in H \mid h_1 \equiv \epsilon\ (\mbox{mod}\ 2) \}. $$
  Then $H_0$ and $H_1$ form a bipartition of $\Cay(H;S-x_1)$. 
Therefore, since $0 \in H_0$,
we know that $d_C(0,h)$ is even for every $h \in H_0$.  Also, since 
$v +  x_1 \in H$ and the
first coordinate of $v + x_1$ is even, we know that $v + x_1 
\in  H_0$.  Hence
$d_C(0,v+ x_1)$ is even, and thus by Proposition~\ref{ham_path}, 
there is a hamiltonian
path in $\Cay \bigl( (\Z_m)^k \bigr)$  from  $0$ to $(v +x_1) - x_1 = v$.
  \end{proof}

Observe that Corollary~\ref{main_even} applies for all~$k$, whenever 
$m$~is even. In
contrast, the following remark shows that if $k = 2$ and $m$~is odd, 
then only $(m+1)/2$ of
the $m$~elements of $H - x_1$ are the terminal endpoint of a 
hamiltonian path that starts
at~$0$.

\begin{rem}
  Assume $k = 2$ and $m$ is odd. Let $v \in H - x_1$, and write $v = 
(v_1,v_2)$ with $0 \le
v_1,v_2 < m$. Then $v_1 + v_2 = m-1$, which is even. Hence $v_1$ and~$v_2$ have 
the same parity. The
work of Curran \cite{CW1} shows that there is a hamiltonian path from 
$0$ to~$v$ in $\Cay
\bigl( (\Z_m)^2 \bigr)$ if and only if $v_1$~is even. Thus, if $m$~is 
odd, then the
assumption that $k \ge 3$ is necessary in Theorem~\ref{main}.
  \end{rem}

We now wish to apply Proposition  \ref{ham_path} in the case that $m$ is odd.  We begin 
by showing that we can find the appropriate hamiltonian cycle.

\begin{lemma} \label{ZaxZb_odd}
  Assume that $m$ is odd and that $n$ is a multiple of~$m$. Let $v \in \Z_m
\times \Z_n$ where $v = (i,j)$ with $0 \le i < m$ and $0 \le j < n$, 
and let $r$ be
the remainder of $i+j$ upon division by~$m$. If either
  \begin{enumerate}
  \item\label{j+r_even} $j+r$ is even, or
  \item\label{jr_konzero} both $j$ and~$r$ are nonzero, or
  \item\label{j_even_konzero} $j$~is even and nonzero,
  \end{enumerate}
  then there is a hamiltonian cycle~$C$ in $\Cay(\Z_m \times \Z_n)$ such that
  $d_C \bigl( 0, v \bigr)$ is even.
  \end{lemma}

\begin{proof}
  As usual, let $x_1 = (1,0)$ and $x_2 = (0,1)$.

  (\ref{j+r_even}) If $j+r$ is even, let $C$ be the hamiltonian cycle
$(x_1^{m-1}, x_2)^n$.
  Then, since
  $$v = ix_1 + jx_2 = j((m-1)x_1 + x_2) + r x_1,$$
  we have that
  $$d_C(0,v) = jm + r \equiv j + r \equiv 0\ (\mbox{mod}\ 2). $$
  Hence $d_C(0,v)$ is even.

  (\ref{jr_konzero}) If $j+r$ is odd, and both $j$ and~$r$ are 
nonzero, let $C$ be the
hamiltonian cycle $(x_2, x_1^{m-1})^n$.  Then, since
  $$v = ix_1 + jx_2
  = (j-1)(x_2 + (m-1)x_1) + x_2 + (r-1)x_1 ,$$
  we have that
  $$d_C(0,v) = (j-1)m + 1 + (r-1) \equiv (j-1) + r \equiv 0\ (\mbox{mod}\ 2). $$
  Hence $d_C(0,v)$ is even.

(\ref{j_even_konzero}) Suppose $j$~is even and nonzero. If $j+r$ is even, then
(\ref{j+r_even}) applies. If not, then $r$ must be odd, so $r \neq 
0$. Therefore
(\ref{jr_konzero}) applies.
  \end{proof}

\begin{prop}\label{m_odd}
  If $m$~is odd and $n \ge 2$, then, for each $v \in (\Z_m)^n$, there is a 
hamiltonian cycle $C$ in
$\Cay \bigl( (\Z_m)^n \bigr)$ such that $d_C(0,v)$ is even.
  \end{prop}

\begin{proof} The proof is by induction on~$n$.

\begin{basecase}
  Assume that $n = 2$.
  \end{basecase}
  Write $v = ix_1 + jx_2$ for some $i$ and $j$ with $ 0 \le i, j \le 
m-1$.  Let $r$ be the
remainder of $i+j$ upon division by~$m$. We may assume that $i+r$ 
and~$j+r$ are both odd,
for otherwise the desired conclusion follows from 
Lemma~\ref{ZaxZb_odd}(\ref{j+r_even})
(perhaps after interchanging $i$ and~$j$). This implies that $i$ 
and~$j$ have the same
parity, so $i+j$ is even. The desired conclusion is obvious if $v = 
0$, so we may assume
that not both of $i$ and~$j$ are~$0$; by symmetry, we may assume $j 
\neq 0$. Now $i + j
\neq 0$ and, because $i+j$ is even, we know $i + j \neq m$; therefore 
$r \neq 0$. Thus, the
desired conclusion follows from Lemma~\ref{ZaxZb_odd}(\ref{jr_konzero}).

\begin{inductionstep}
  Let $n \ge 3$. Assume for every vertex~$w$ of $\Cay \bigl( (\Z_m)^{n-1} \bigr)$, 
there is a hamiltonian cycle $C$ such that $d_C(0,w)$ is even.
  \end{inductionstep}
  Write $v = (v_1,\ldots,v_n)$. We may assume $v \neq 0$, for 
otherwise the desired
conclusion is obvious. By symmetry, we may assume that $v_n \neq 0$. Let 
$w = (v_2,\ldots,v_n)$.
By induction, there is a hamiltonian cycle~$C_0$ in $\Cay \bigl( 
(\Z_m)^{n-1} \bigr)$, such
that $d_{C_0}(0,w)$ is even. Let $j = d_{C_0}(0,w)$. Then there is an 
embedding $\phi$ of
$\Cay(\Z_m \times \Z_{m^{n-1}})$ as a spanning subdigraph of $\Cay 
\bigl( (\Z_m)^n \bigr)$,
such that $\phi(0) = 0$ and $\phi(v_1,j) = v$. Because $j$ is even and nonzero,
Lemma~\ref{ZaxZb_odd}(\ref{j_even_konzero}) implies that there is a 
hamiltonian cycle~$C$ in
$\Cay(\Z_m \times \Z_{m^{n-1}})$, such that $d_C \bigl( 0, (v_1,j) 
\bigr)$ is even. Then
$\phi(C)$ is a hamiltonian cycle in  $\Cay \bigl( (\Z_m)^n \bigr)$, such that
$d_{\phi(C)}(0,v)$ is even.
  \end{proof}

Combining the results of Propositions \ref{ham_path} and \ref{m_odd} together with 
Corollary \ref{main_even}, we now give the proof of our main result.

\begin{proof}[Proof of Theorem~\ref{main}]
  If $m$ is even, then the desired conclusion follows by 
Corollary~\ref{main_even}.
Thus, we may assume that $m$ is odd. Since $\Cay(H;S-x_1)$ is 
isomorphic to $\Cay
\bigl( (\Z_m)^{k-1} \bigr)$, Proposition~\ref{m_odd} implies that 
there is a hamiltonian
cycle $C$ in $\Cay(H;S-x_1)$ such that $d_C(0,v)$ is even. The 
desired conclusion now
follows from Proposition~\ref{ham_path}.
  \end{proof}

\begin{ack}
  D.W.\ was partially supported by a grant from the National Science Foundation
(DMS--9801136).
  \end{ack}

\end{document}